\documentclass{amsproc}
\usepackage{amscd}
\usepackage{eucal}
\newcommand{\R}{{\mathbb R}}
\newcommand{\C}{{\mathbb C}}
\newcommand{\h}{{\rm H}}
\newcommand{\Hol}{{\rm Hol}}
\newcommand{\Hom}{{\rm Hom}}
\newcommand{\CP}{{\mathbb C}P}
\newcommand{\M}{{\overline{\mathcal M}}}
\newcommand{\pt}{{\rm pt}}
\newcommand{\Q}{{\mathbb Q}}
\newcommand{\T}{{\mathbb T}}
\newcommand{\Tr}{{\rm Tr}}
\newcommand{\Z}{{\mathbb Z}}
\newcommand{\D}{{\bf D}}
\newcommand{\e}{{\bf e}}
\newcommand{\s}{{\bf S}}
\newcommand{\MU}{{\bf MU}}
\newcommand{\half}{\textstyle{\frac{1}{2}}}
\newcommand{\eighth}{\textstyle{\frac{1}{8}}}
\begin{document}

\title{Quantum generalized cohomology}

\author{Jack Morava}
\address{Department of Mathematics, Johns Hopkins University, 
Baltimore, Maryland 21218}
\email{jack@math.jhu.edu}
\thanks{Author was supported in part by NSF Grant \#9504234.}
\date{15 July 1998}
\begin{abstract} We construct a ring structure on complex 
cobordism tensored with $\Q$, which is related to 
the usual ring structure as quantum cohomology is related to ordinary 
cohomology. The resulting object defines a generalized 
two-dimensional topological field theory taking values in a
category of spectra.\end{abstract}

\maketitle

\section*{Introduction}

The conclusion of this paper is that the theory of 
two-dimensional topological gravity has a remarkably
straightforward homotopy-theoretic interpretation in terms of a
generalized cohomology theory, completely analogous to the more
familiar interpretation of quantum ordinary cohomology as a
topological field theory. Two-dimensional gravity originated in
attempts to integrate over the space of metrics on a Riemann
surface; it was reformulated by Witten in terms of an
algebra of generalized Miller-Morita-Mumford characteristic
classes for surface bundles. In the interpretation proposed here,
this algebra is the coefficient ring of a generalized cohomology
theory, and topological gravity becomes a
topological-field-theory-like functor, which assigns invariants
to families of algebraic curves just as a classical topological
field theory assigns invariants to individual curves; in this it
resembles algebraic $K$-theory, which assigns homotopy-theoretic
invariants to families of modules over a ring. 

Here is an outline of the argument. After a preliminary
section which collects some background information, we define a
generalized topological field theory in \S 2 in terms of
(homotopy classes of) maps $$\tau_{g}^{n} : \M_{g}^{n}
\rightarrow {\bf M}^{\wedge n}$$ from a compactified moduli space
of curves marked with $n$ points, to $n$-fold powers of a 
module-spectrum {\bf M}. These maps preserve a monoidal structure
(defined geometrically in the domain by glueing curves together
at marked points, but defined algebraically in the range); in the
language of [14 \S 3, 17 \S 1.7], $\tau_{*}^{*}$ is a representation
of a certain cyclic operad. The existence of such a
representation entails the existence of a (quantum)
multiplication on the module-spectrum {\bf M} (cf.\ \S 2.4); in
familiar cases this is the multiplicative structure defined by
the WDVV equation, and when $n=0$ we recover Witten's
tau-function for the moduli space of curves.

That topological gravity and quantum cohomology are
closely related is clear from [37], but I suspect that the
simplicity of the underlying geometry is not widely understood.
The main technical lemma (\S 2.1) is a kind of splitting theorem
(cf.\ [21, 33]) for generalized Gromov-Witten classes; it is quite
natural, in light of Kontsevich's ideas about stacks of stable maps. 
These moduli objects have not
yet been shown to be smooth for curves of all genus for any
variety more complicated than a point [21 \S 1], so our main
example is still conjectural, but in fact smoothness is
much more than we need; the constructions of this paper require
only that the generalized Gromov-Witten maps of \S 1.3 be local complete
intersection morphisms. This is a convenient working hypothesis,
which can be weakened further by elaborating the cohomological
formalism; but this paper concerned with the consequences of this
assumption, not with proving it. For curves of genus zero 
the moduli stacks are known to be smooth, if the target manifolds 
are convex in a suitable sense [4]; this leads to a simple 
proof (\S 2.2) of the associativity of the quantum multiplication,
and when the defining variety is a point, we can calculate the 
corresponding coupling constant (\S 2.3).  

Two short appendices discusses some related
issues. In particular, there is reason to think that
the Virasoro algebra is a ring of `quantum generalized cohomology
operations' for the main example. To state this more
precisely requires a short digression about representions of the
group of antiperiodic loops on the circle, which is included as
the first appendix. A second appendix, added after the publication 
of this paper [in the Proceedings of the Hartford-Luminy Conference
on the Operad Renaissance, Contemporary Math. 202 (1997) 407-419],  
outlines a slight generalization of the main construction of this
paper in terms of the physicists' `large phase space' of deformations
of quantum cohomology. 

It is a pleasure to thank Ralph Cohen, Yuri Manin, and
Alexander Voronov for helpful conversations about the content of
this paper; that the paper exists at all, however, is the
consequence of helpful conversations with Graeme Segal and Edward
Witten. 

\section{Notation and conventions}

\subsection{} Let $V$ be a simply-connected projective
smooth complex algebraic variety of real dimension $2d$, with
first Chern class $c_{1}(V)$, and let $\h$ denote its second
integral homology group $H_{2}(V,\Z)$. We will use a rational
version $$\Lambda = \Q[\h \times \Z]$$ of the Novikov ring [25 \S
1.8] of $V$: its elements are Laurent polynomials $$\sum_{k \in
\Z, \alpha \in \h} c_{\alpha,k} \alpha \otimes v^{k}$$ with
coefficients $c_{\alpha,k} \in \Q$. This ring has a useful
grading, in which $v$ has (cohomological) degree two, and
$\alpha$ has degree $2\langle c_{1}(V), \alpha \rangle$. We will
also use the notation $$v_{(k)}  = \frac {v^{k}}{k!}$$ for the
$k$-th divided power of $v$. If $u: \Sigma \rightarrow  V$ is a
map from a connected oriented surface to $V$, then the degree of
$u$ is the class $u_{*}[\Sigma] \in \h$, where $[\Sigma] \in
H_{2}(\Sigma,\Z)$ is the fundamental class of the surface.

\subsection{} 
$MU_{*}(X)$ will denote the complex bordism of a
$CW$-space $X$, and $\MU_{\Lambda}$ is the spectrum representing
the homology theory defined on base-pointed finite complexes by
$$X \mapsto MU_{*}(X) \hat \otimes \Lambda = [S^{*},X \wedge
\MU_{\Lambda}];$$ the tensor product has been completed, so
$MU^{*}(\pt) \hat \otimes \Lambda$ is the graded ring of formal
Laurent series in $v$, with coefficients from the graded group
ring $MU^{*}(\pt)[\h]$. It will be convenient to write
$\MU_{\Lambda}(V)$ for the function spectrum
$F(V^{+},\MU_{\Lambda})$. [The superscript $+$ indicates the
addition of a disjoint basepoint, but this refinement will often 
be omitted when the space is already encumbered with superscripts.] 
The fiber product of spaces (or
schemes) $X$ and $Y$ over $Z$ will be denoted $X \times_{Z}Y$,
and the product of $MU^{*}_{\Lambda}(X)$ and
$MU^{*}_{\Lambda}(Y)$ over $MU^{*}_{\Lambda}(Z)$ will be denoted
$MU^{*}_{\Lambda}(X) \otimes_{Z} MU^{*}_{\Lambda}(Y)$. The
spectrum $\MU_{\Lambda}(V)$ has an $\MU_{\Lambda}$-algebra
structure, and $$\MU_{\Lambda}(V^{n}) = \MU_{\Lambda}(V)
\wedge_{\MU_{\Lambda}} \dots \wedge_{\MU_{\Lambda}}
\MU_{\Lambda}(V)$$ is its $n$-fold Robinson smash power [31] over
$\MU_{\Lambda}$. There is a map $$\Tr_{V} : \MU_{\Lambda}(V)
\rightarrow \MU_{\Lambda}$$ of $\MU_{\Lambda}$-module
spectra, which represents the transfer map $$MU_{\Lambda}^{*}(X
\wedge V^{+})\rightarrow MU_{\Lambda}^{*-2d}(X)$$ defined by the
(complex oriented) projection $V \rightarrow  \pt$, followed by 
multiplication with $v^{d}$ (to shift dimensions).

According to Quillen, a proper complex-oriented map
$\Phi : P \rightarrow  M$ between smooth manifolds defines an
element $[\Phi]$ of $MU^{k}(M)$, where $k$ is the codimension of
$\Phi$; more generally, a suitably oriented map between geometric
cycles, with enough of a normal bundle to possess rational Chern 
classes, will define an element of $MU^{k}_{\Q}(M)$.
Contravariant maps in cobordism are defined by fiber
products, while covariant map are defined by the obvious
compositions. Finally, the bilinear form $$b_{V} :
\MU_{\Lambda}(V) \wedge \MU_{\Lambda}(V) \rightarrow
\MU_{\Lambda}(V) \rightarrow \MU_{\Lambda} $$ is the composition
of the trace with the multiplication map of $\MU_{\Lambda}(V)$.

The graded ring $MU^{*}_{\Lambda}(V)$ is a technical
replacement for $MU^{*}_{\Q[v,v^{-1}]}(V \times \h)$, which is in
some ways more natural; but the latter ring does not help with
the usual convergence problems, which (in the present framework)
are consequences of the failure of the map $\h \times \h
\rightarrow \pt$ to be proper. 

\subsection{}
A (marked) algebraic curve is stable if its group of
automorphisms is finite; $\M_{g}^{n}$ will denote the
Deligne-Mumford-Knudsen space of such curves of arithmetic genus
$g$, marked with $n$ ordered smooth points. These spaces are
compact orbifolds, of complex dimension $3(g-1) + n \geq 0$;
cases of low genus are thus sometimes exceptional. It
is useful to understand $n$ to be a finite ordered set (or
ordinal number), so that permutations of $n$ can act on $\M_{g}^{n}$.  
More generally, $\M_{g}^{n}(V,\alpha)$ will
denote the stack [4 \S 3, 12] of stable maps of degree $\alpha$
from a curve of genus $g$ marked with $n$ ordered smooth points,
to $V$; there is a morphism from $\M_{g}^{n}(V,\alpha)$ to
$\M_{g}^{n}$ which assigns to a map (the stabilization of) its domain, 
and there is a morphism to $V^{n}$ which evaluates a map at the marked
points. The product of these is a perfect (finite Tor-dimension 
[11 II \S 1.2]) proper morphism $$\Phi_{V,g,\alpha}^{n} : \M_{g}^{n}(V,\alpha)
\rightarrow \M_{g}^{n} \times V^{n} $$ of stacks. At a point
$u:\Sigma \rightarrow V$ of $\M_{g}^{n}(V,\alpha)$ defined
by a smooth $\Sigma$, the relative tangent space to $\M_{g}^{n}(V,\alpha)$ 
over $\M_{g}^{n}$ defines a $K$-theory
class $$[H^{0}(\Sigma,u^{*}T_{V})] - [H^{1}(\Sigma,u^{*}T_{V})]$$
of complex dimension $d(1-g) + \langle c_{1}(V),\alpha \rangle$, where
$d$ is the complex dimension of $V$. It
seems likely that under reasonable hypotheses this map will be
a local complete intersection morphism of stacks, and in particular
that the $K$-theory class of its cotangent complex at a singular point
$(\Sigma, u)$ will equal the holomorphic Euler class of the pullback of
$u^{*}T_{V}$ to the normalization of $\Sigma$. [In fact, Kontsevich
[19 \S 1.4] has already sketched something very close to a local complete
intersection structure for this morphism.].
Similarly, let  $\epsilon_{V,g,\alpha}^{n}(k)$ denote the proper
complex-oriented map from $\M_{g}^{n+k}(V,\alpha)$ to
$\M_{g}^{n}(V,\alpha)$ which forgets the final $k$ marked
points, and define $$ \Phi_{V,g,\alpha}^{n}(k) =
\Phi_{V,g,\alpha}^{n} \circ \epsilon_{V,g,\alpha}^{n}(k):
\M_{g}^{n+k}(V,\alpha) \rightarrow \M_{g}^{n} \times V^{n} \;.$$ 
There are also generalizations $$\mu_{V}^{s} : \M_{g}^{r+s}(V,\alpha) 
\times_{V^{s}} \M_{h}^{s+t}(V,\beta) \rightarrow \M_{g+h+s-1}^{r+t}(V,\alpha
+ \beta) \;,$$ of Knudsen's glueing morphisms [17]. All these 
maps represent natural transformations between moduli functors,
so their normal bundles are reasonably accessible. In diagrams
below, complicated subscripts and superscripts will be supressed 
when they are redundant in context.

\section{Generalized topological field theories}

\subsection{} We will be interested in generalized Gromov-Witten
invariants defined by the cobordism classes of these morphisms. {\bf I 
will assume that these maps are local complete
intersection morphisms}; such
maps between (possibly singular) varieties have most of the topological
transversality properties of maps between smooth manifolds. In 
particular, they have well-behaved Gysin homomorphisms and normal bundles 
[1 IV \S 4] and they thus define elements in complex cobordism tensored with 
the rationals. [A variant approach is discussed below in \S 2.5.] 
Our generalized Gromov-Witten invariants are the classes
$$\phi_{V,g}^{n}(k) = \sum_{\alpha \in H}
[\Phi_{V,g,\alpha}^{n}(k)] \alpha \otimes v_{(k)} \in
MU_{\Lambda}^{2d(n+g-1)}(\M_{g}^{n} \times V^{n}) \; ;$$
permutations of $k$ define cobordant elements. Summing these
classes over $k$ defines $$\tau_{V,g}^{n} = v^{-d(n+g-1)} \sum_{k
\geq 0} \phi_{V,g}^{n}(k) \in MU_{\Lambda}^{0}(\M_{g}^{n}
\times V^{n}) \;;$$ the convergence problems mentioned in the
preceding section do not appear when the function $\alpha \mapsto
[\Phi_{V,g,\alpha}^{n}(k)]$ is supported in a proper cone in $H$.
The tau-function $\tau_{V,g}^{n}$ can be interpreted
geometrically as the cobordism class of the `grand canonical 
ensemble' of maps from a curve of genus $g$ marked with $n$
ordered smooth points, together with an indeterminate number of
further distinct smooth (unordered) points, to $V$, but we will
be more concerned with the homotopy class $$\tau_{V,g}^{n} : 
\M_{g}^{n} \rightarrow \MU_{\Lambda}(V^{n})$$ it defines.\medskip

\noindent {\sc Proposition} 2.1: The diagram
 $$\begin{CD}\M_{g}^{r+s} \wedge \M_{h}^{s+t} &@>\mu_{s}>> & \M_{g+h+s-1}^{r+t}
\\@VV{\tau \wedge \tau}V && @VV {\tau}V \\ \MU_{\Lambda}(V^{r+s})
\wedge \MU_{\Lambda}(V^{s+t}) &@>b_{V}^{s}>>& \MU_{\Lambda}
(V^{r+t}) \end{CD}$$ commutes up to homotopy; alternately,
$$\mu^{*}_{s}(\tau_{V,g+h+s-1}^{r+t}) = v^{sd} \Tr_{V}^{s}
(\tau_{V,g}^{r+s} \otimes_{V^{s}} \tau_{V,h}^{s+t}) \;.$$ 
{\it Sketch Proof}, under the standing hypothesis
above: The general case reduces by
induction to the case $s=1$, which can be stated as a coproduct
formula $$\mu^{*}\phi_{V,g+h}^{r+t}(k) = v^{d} \sum_{i+j=k}
\Tr_{V}(\phi_{V,g}^{r+1}(i) \otimes_{V} \phi_{V,h}^{1+t}(j)) .$$
This can be reformulated as the assertion that the two diagrams 
$$\begin{CD}\bigsqcup_{i+j=k} \M_{g}^{i+r+1} (V,\alpha)
\times_{V} \M_{h}^{1+t+j}(V,\beta)& @>\mu_{V}>>& \M_{g+h}^{r+t+k}(V,\alpha + 
\beta) \\@VV{ \bigsqcup \epsilon(i) \times \epsilon(j)}V& &@VV{\epsilon(k)}
V\\ \M_{g}^{r+1}(V,\alpha) \times_{V} \M_{h}^{1+t}(V,\beta)
&@>\mu>>&\M_{g+h}^{r+t}(V,\alpha + \beta)\end{CD}$$ and
$$\begin{CD}\bigsqcup_{\alpha + \beta = \gamma} \M_{g}^{r+1}(V,\alpha) 
\times_{V} \M_{h}^{1+t}(V,\beta)& @>\mu_{V}>>&\M_{g+h}^{r+t}(V,\gamma)\\@VV 
{\bigsqcup \Phi_{\alpha} \times \Phi_{\beta}}V & &@VV{ \Phi_{\gamma}}V
\\ \M_{g}^{r+1} \times V^{r+1+t} \times \M_{h}^{1+t}
&@>\mu \times \Tr_{V}>>&\M_{g+h}^{r+t} \times V^{r+t}
\end{CD}$$ are fiber products; the claim follows by stacking the
first diagram on top of the second. The bottom diagram describes
the stable maps of decomposable curves in terms of the
restrictions to their components. In the top diagram, the union
is to be taken over partitions of the set with $k$ elements into
subsets of cardinality $i$ and $j$; on the level of functors,
this diagram asserts that the ways of sprinkling points on a
curve decomposed into two components correspond to the ways of
sprinkling points on the components separately.

The last assertion is not entirely straightforward,
because forgetting marked points may destabilize a genus zero
component of a stable marked curve; the morphism $\Phi_{V}(k)$
will blow such components down to points. The union in the upper
left corner of the diagram is thus not necessarily disjoint: the
fiber product is obtained from the disjoint union by
identification along certain divisors [28 \S 3]. The point,
however, is that the cobordism classes are defined by maps rather
than by subobjects; the fiber product class is equivalent to the
sum of the classes defining the disjoint union.

\subsection{} Proposition 2.1 states that the triple 
$(\tau_{V*}^{*}, \MU_{\Lambda}(V), b_{V})$
defines a topological field theory which takes values in the
category of $\MU_{\Lambda}$-module spectra, where the usual
monoidal structure defined by tensor product of modules over a
ring is replaced by the smash product of module-spectra
over a ring-spectrum. The domain of this generalized topological
field theory is the monoidal category {\it (Stable Curves)}
with finite ordered sets as objects; morphisms are finite unions
of marked curves. [This category, however, does not possess
identity maps for its objects.] Both the domain and range of the
generalized topological field theory are topological categories,
and $\tau_{V*}^{*}$ defines a homotopy class of maps from the
space of morphisms of the domain category to the space of
morphisms of the range. These homotopy classes preserve the
composition of morphisms and thus define a functor. In this
generality, we need a version of Proposition 2.1 for Knudsen
glueing of two points on a connected curve, but the changes required
for this are minor.

The construction involves the bilinear map $b_{V}$,
but it has not otherwise used the multiplicative structure on
$\MU_{\Lambda}(V)$. In fact the morphism $$\tau_{V,0}^{3} : S^{0}
= \overline M_{0}^{3} \rightarrow \MU_{\Lambda}(V^{3}) $$
defines a composition 
$$\begin{CD} *_{V} : \MU_{\Lambda}(V) \wedge S^{0} \wedge
\MU_{\Lambda}@>id \wedge \tau_{0}^{3} \wedge id  >>
\MU_{\Lambda}(V^{5}) @> b_{V}^{\otimes 2} \wedge id>>
\MU_{\Lambda}(V) \;, \end{CD}$$ and Proposition 2.1 has the
following \bigskip

\noindent {\sc Corollary} 2.2 : The pair
$(\MU_{\Lambda}(V),*_{V})$ is a homotopy commutative and
homotopy associative ring-spectrum. \bigskip

\noindent{\it Sketch proof}: Smashing
the morphism $$ A(2,2) := (id \wedge b_{V} \wedge id)(\tau_{0}^{2+1}
\wedge \tau_{0}^{1+2}) : S^{0} = \M_{0}^{2+1} \wedge
\M_{0}^{1+2} \rightarrow \MU_{\Lambda}(V^{4}) $$ with
the identity map of $\MU_{\Lambda}(V^{3})$ defines a map from 
$\MU_{\Lambda}(V^{3})$ to $\MU_{\Lambda}(V^{7})$; arranging the
seven copies of $V$ into pairs and applying the trace map $b_{\Lambda}$ 
three times defines a collection of maps from $\MU_{\Lambda}(V^{3})$ to
$\MU_{\Lambda}(V)$ indexed by the possible groupings of the factors. By
our conventions $2+1$ and $1+2$ are isomorphic but not equal, so the 
notation for this associator class
emphasizes that it depends on four points partitioned
into two subsets, each containing two items. Ignoring obvious
involutions, there are three different partitions, corresponding
to the maps $\pi_{0}$, $\pi_{1}$, $\pi_{\infty}$ from the 0-manifold
$\M_{0}^{2+1} \times \M_{0}^{1+2}$ to $\M_{0}^{4}$ which send it to a 
degenerate curve of genus zero with two irreducible components, each 
carrying two marked points
(aside from the node); the cross-ratio identifies these
configurations with the standard points 0,1, and $\infty$ on the projective
line. To verify associativity it
suffices to show that the homotopy class $A(2,2)$ is independent
of the way the four points are partitioned into pairs; but by the
proposition, $A(2,2) = \tau_{0}^{4} \circ \pi_{i}$ factors through
$\overline M_{0}^{4}$, where the three maps $\pi_{i}$ become
homotopic. 
 
\subsection{} The morphism $$\tau_{0}^{3} = v^{-2d} \sum_{k
\geq 0} [\M_{0}^{k+3}(V) \rightarrow V^{3}]v_{(k)}$$
defining this quantum multiplication is essentially the
Gromov-Witten potential [21].  Because the moduli spaces
$\M_{g}^{n}$ are not defined when $3(g-1)+n$ is negative,
however, it is not clear that the resulting multiplicative
structure on $\MU_{\Lambda}(V)$ possesses a unit. The class
$$q_{V} := 1*_{V}1 = v^{-2d}(b_{V} \otimes id)(\tau_{V,0}^{3})$$
is the coupling constant for the topological field theory defined
by $MU_{\Lambda}^{*}(V)$; this theory assigns to a connected
surface of genus $g$ with one boundary component, the $2g$th
power of 1 with respect to the product $*_{V}$. \bigskip

\noindent
{\sc Corollary} 2.3 : When $V$ is a point, the 
coupling constant of the resulting topological field theory is
$$q = \sum_{k \geq 0} [\M_{0}^{k+3}]v_{(k)}  \in
MU^{0}_{\Lambda}(pt) ,$$ and the quantum product in
$MU_{\Lambda}^{*}(pt)$ is $x*y = qxy$. \bigskip

In this formula $[\M_{0}^{k+3}]$ is the cobordism
class of the manifold of configurations of $k+3$ points on a
curve of genus zero; the sum is thus a cobordism analogue of
Manin's Hodge-theoretic invariant $\varphi$ [26 \S 0.3.1]. With this
structure, $(MU_{\Lambda}^{*}(\pt),*_{\pt})$ is isomorphic to
$MU_{\Lambda}^{*}(\pt)$ with its usual multiplication, by a
homomorphism which sends $1^{*2g}$ to $q^{2g}$. More generally,
the operation $x \mapsto 1*_{V}x$ is a module isomorphism, and it
seems reasonable to hope that $(1*_{V})^{-1}(1)$ will be a unit
for $*_{V}$.

Knudsen glueing defines a pair-of-pants product
$$\mu^{+}: \M_{g}^{1} \wedge \M_{h}^{1} \rightarrow \M_{g}^{1} \wedge 
\M_{0}^{3} \wedge \M_{h}^{1} \rightarrow \M_{g+h}^{1}$$ and it follows 
from \S 2.1 and the arguments above that the diagram 
$$\begin{CD} \M_{g}^{1} \wedge \M_{h}^{1} @>\mu^{+}>> \M_{g+h}^{1} \\ @VV
{ \tau \wedge \tau}V@VV{ \tau}V \\ \MU_{\Lambda} \wedge \MU_{\Lambda} @>*_{V}>>
\MU_{\Lambda} \end{CD}$$ is homotopy-commutative; in other words,
$$\tau_{V*}^{1} : \M_{*}^{1} \rightarrow \MU_{\Lambda}$$
is a kind of homomorphism of monoids. This is probably the most
intuitive way to think of the product in quantum cohomology, but
from the present point of view it is a conclusion, not a
definition.
          
\subsection{} A generalized topological field theory has an
associated theory of topological gravity, which assigns
invariants to proper flat families of stable curves. Such a
family $Z$, say of topological type $(g,n)$, is defined by its
classifying map to $\M_{g}^{n}$; the pullback of
$\tau_{V,g}^{n}$ along this morphism defines a class $\tau_{V}(Z)
\in [Z^{+},\MU_{\Lambda}(V^{n})]$. If (for simplicity) we assume
that $V$ is a point, and write $[Z] \in H_{*}(Z)$ for the
fundamental class of $Z$, then the image $\tau_{*}(Z)$ of $[Z]$
in $H_{*}(\MU,\Lambda)$ under the map induced on homology by
$\tau(Z)$ is a kind of absolute invariant of the family, obtained
by integrating $\tau(Z)$ over $[Z]$. In particular, the vacuum 
morphism $$0 \rightarrow 0$$ is defined by the family of
arbitrary finite unions of unmarked stable curves; the infinite
symmetric product ${\rm SP}^{\infty}(\bigsqcup_{g \geq 0} \M_{g}^{0})$ 
is a rational model for its parameter space. The
resulting absolute invariant $$\tau = \exp (\sum_{g \geq 0}
\tau_{*}(\M_{g}^{0})) \in H_{*}(\MU,\Q[v,v^{-1}])$$ is
Witten's tau-function for two-dimensional topological gravity
[24]. The point is that the characteristic number homomorphism
$$MU^{*}(M) \rightarrow H^{*}(M,H_{*}(\MU))$$ sends the class
$[\Phi:P \rightarrow M]$ to a sum of the form $\sum_{I}
\Phi_{*}c_{I}(\nu) t^{I}$, where $\nu$ is the stable normal
bundle of $\Phi$, $c_{I}(\nu)$ is a certain polynomial (indexed
by $I = i_{1},\dots$) in its Chern classes, $t^{I} =
t_{1}^{i_{1}} \dots$ is a product of elements in
$H_{*}(\MU) = \Z [t_{i}|i \geq 1]$, and $\Phi_{*}$ is the
covariant (Gysin) homomorphism induced by $\Phi$. The stable
normal bundle of $\Phi(k)$ is inverse to the tangent bundle along
the fibers, which is the sum of the $k$ line bundles defined by
the tangent space to the universal curve at its $k$ marked
points, so its Chern polynomials can be expressed as polynomials
in the Chern classes of these line bundles. Under the pushdown
$\Phi(k)_{*}$, these become polynomials in the Mumford classes in
the cohomology of $\M_{g}^{0}$.

\subsection{} Some aspects of Kontsevich-Witten theory suggest that 
Gromov-Witten invariants can be defined more naturally in $K$-theory
than in ordinary cohomology. The algebraic $K$-theory of a
reasonable stack, tensored with the rationals, agrees with the
algebraic $K$-theory of its quotient space [15 \S 7], but (perhaps
because of this) the $K$-theory of stacks seems to have received
little attention otherwise. The following assumes that the
standard direct image construction for perfect proper maps of
schemes [36 \S 3.16.4] generalizes to stacks.
 
Let $\pi_{g,n} : C_{g}^{n} \rightarrow \M_{g}^{n}$ be the
universal stable curve; because the range is smooth, this is a
perfect proper morphism.  Let $$C_{g,\alpha}^{n}(V) := \M_{g,\alpha}^{n}(V) 
\times_{\M_{g}^{n}} C_{g}^{n}\; ; $$ from now
on I will supress the subscripts. The projection $\bar \pi : C(V)
\rightarrow \M(V)$ to the first factor,
being the pullback of a perfect morphism, is again perfect. Let $U
: C(V) \rightarrow V$ be the universal evaluation morphism; the
vector bundle $U^{*}T_{V}$ defines an element of $K(C(V))$, and
its hypothetical direct image $\bar \pi_{*} U^{*}T_{V} :=
\nu(\Phi_{V}) \in K(\M(V))$ is a reasonable candidate
for the normal bundle of $\Phi_{V}$.
 
Because $\Phi_{V}$ is itself proper and perfect, we can
define generalized Gromov-Witten classes $$\sum_{I}
\Phi_{V*}m_{I}(\nu(\Phi_{V}))t^{I} \in K^{*}(\M_{g}^{n}
\times V^{n}) \otimes_{K}K_{*}\MU \;, $$ where $I$ is a multiindex
as above, $t^{I} = \prod_{k} t_{k}^{i_{k}}$ is a basis for
$K_{*}\MU$, and $m_{I}$ denotes the $K$-theory characteristic
class associated to the monomial symmetric function by the
correspondence which assigns gamma operations [5 V \S 3] to the elementary
symmetric functions. By the Hattori-Stong theorem, such a sum can
be identified with a class in the localization $MU^{*}(\M_{g}^{n} \times 
V^{n})[\CP(1)^{-1}]$ of complex cobordism.
 
If $\Phi_{V}$ is a local complete intersection morphism, this
approach to defining Gromov-Witten invariants agrees with the
definition in \S 2.1. In any case some such hypothesis seems to be needed
to make the arguments of Prop.\ 2.1 work.

\section{Some questions}
\subsection{}
Kontsevich and Witten [20, 24, 37] show that the
tau-function for the vacuum state of two-dimensional topological 
gravity is a lowest weight vector for a certain
representation of the Virasoro algebra. This Lie algebra
bears a striking resemblance to the Lie algebra defining the 
Landweber-Novikov algebra of operations in complex cobordism, but
the relation between these two structures is not well-understood.
I have included as an appendix a construction for the 
Kontsevich-Witten representation, starting from a representation
of a certain loop group of antiperiodic functions on the circle,
following [8]. One point of the appendix is 
that the representation theory of this loop group is essentially
trivial. 

On the other hand, the usual complex cobordism functor
takes values in the monoidal category of $\Z/2\Z$-graded 
{\bf G}-equivariant sheaves over the moduli scheme Spec
$MU^{*}(\pt)$ of formal group laws, with the Landweber-Novikov
group {\bf G} of formal coordinate transformations acting by
change of coordinate; but this category is equivalent, after
tensoring with $\Q[v,v^{-1}]$, to the category of $\Z/2\Z$-graded
vector spaces. It therefore seems not completely unreasonable to
conjecture that the group of antiperiodic loops is a kind of
motivic group for the generalized quantum cohomology theory
defined by $\MU_{\Lambda}$. 

As for products in $V$, the functor $\tau_{V*}^{*}$ 
seems to behave very naturally (cf.\ [19]); in particular, it 
is reasonable to expect that $$\tau_{V_{0} \times V_{1},g}^{n} =
\tau_{V_{0},g}^{n} \otimes_{MU^{*}_{\Lambda}} \tau_{V_{1},g}^{n}\;.$$

\subsection{} The work in this paper was originally motivated
by a desire to understand topological gravity and quantum
cohomology from the point of view of Floer homotopy theory [6,7],
but such questions have been supressed here. It may be helpful,
however, to observe that the circle group $\T$ acts on the
universal cover $\widetilde{LV}$ of the free loopspace of $V$,
with $V \times \h$ as fixed point set, so we can think of
$MU_{\Lambda}^{*}(V)$ as its $t_{\T}MU_{\Q}^{*}$-cohomology
[16 \S 15]. The Tate cohomology
$t_{\T}MU_{\Q}^{*}(\widetilde{LV})$ is a rough approximation to the 
Floer $\MU$-homotopy type of $\widetilde{LV}$, and we might hope 
to understand the relation between these invariants as a 
localization theorem for Tate cohomology. 

More specifically, given a compact pointed Riemann surface $(\Sigma,x)$,
let $$(D,0) \rightarrow (\Sigma,x)$$ be a holomorphically embedded closed
disk; the boundary $\partial D$ separates the surface into components $\bar
\Sigma_{0}$ and $\bar \Sigma_{\infty}$, with $x$ the point at infinity, as in
[30 \S 8.11]. Let $\Hol(\Sigma,D;V)$ denote the space of continuous
maps from $\Sigma$ to $V$ which are holomorphic on $\Sigma_{0}$ and
$\Sigma_{\infty}$. This is a manifold, with tangent space $$H^{0}
(\bar \Sigma_{0},u^{*}_{0}T_{V}) \oplus H^{0}(\bar \Sigma_{\infty},
u^{*}_{\infty}T_{V})$$ at $u \in \Hol(\Sigma,D;V)$; here
the sections of the pullback bundles are to be holomorphic on the interior
and smooth on the boundary. Restriction to the boundary defines a map $u
\mapsto \partial u$ to the free loopspace of $V$, but the homotopy class of
$u_{\infty}$ defines a canonical contraction of $\partial u$, so this
restriction map factors naturally through a lift to the universal cover
$\widetilde {LV}$ of $LV$.
 
This map is Fredholm, with index equal to the holomorphic
Euler characteristic of $u^{*}T_{V}$. [More precisely: since $u$ will usually
not be holomorphic, $u^{*}T_{V}$ can't be expected to be holomorphic either;
but $u^{*}T_{V}$ restricts to a holomorphically trivial bundle on an annulus
containing $\partial D$, so $u^{*}_{0}T_{V}$ extends to a holomorphic bundle
$\tilde u^{*}T_{V}$ on $\Sigma$. Then $\chi(\tilde u^{*}T_{V})$ is the
index at $u$.] Moreover, away from maps which collapse $\partial D$ to a
point, this map appears to have a good chance to be proper.
 
We can elaborate this construction, by considering the space
$\Hol(\Sigma_{\hat x},V)$ of holomorphic disks in $\Sigma$
centered at $x$, together with a map to V, continuous and holomorphic away
from $\partial D$ as above; since we're enlarging things, we may as well
include trivial disks too. This thickening has the same homotopy type as the
preceding space, but now $\T$ acts by rotating loops. More generally, we
can allow the moduli of $\Sigma$ to vary as well, thus defining a space of
maps over a thickening of the moduli space $\M$. Restriction maps
this space to a similar thickening of $\M \times \widetilde {LV}$,
defining a candidate for a proper
$\T$-equivariant Fredholm map, and thus an element of
$MU^{-2\chi}_{\T}(\M \times \widetilde {LV})$, which restricts to the classical
Gromov-Witten invariant at the fixed point set of $\T$.

\section*{Appendix I : $MU^{*}(\pt)$ as a 
Virasoro-Landweber-Novikov bimodule}

The Virasoro algebra is the Lie algebra of a central
extension of the group $\D$ of diffeomorphisms of a circle; it is
true generally [30 \S 13.4] that $\D$ acts projectively on a
positive energy representation of a loop group, and in this
appendix I will sketch the construction of an action of the
double cover $\D(2)$ of $\D$ on the basic representation of the
twisted loop group $$L\T\mbox{wist} = \{f \in L\T | f = \iota(f) 
\}$$ of functions from the circle $\R/\Z$ to $\T = \{ z \in \C |
|z| = 1\}$ which are invariant under the involution $\iota(f)(x)
= f(x + \half)^{-1}$.

Because the loop functor preserves fibrations, the
exact sequence of the exponential function $\e(x) = e^{2\pi i x}$
yields an exact sequence $$0 \rightarrow \Z \rightarrow 
L\R \rightarrow L\T_{0} \rightarrow 0$$ of abelian groups with
involution, the group on the right being the identity component
of the group of untwisted loops. The associated exact sequence
$$0 \rightarrow L\R^{\Z/2\Z} \rightarrow L\T\mbox{wist}
\rightarrow H^{1}(\Z/2\Z,\Z) \rightarrow 0$$ of cohomology groups
presents the antiperiodic loops as a canonically split extension
of the group $\Z/2\Z$ of constant loops with value plus or minus
one, by a vector space of antiperiodic functions.

Now $L\T_{0}$ contains a subgroup $\T$ of constant
loops, and $\D$ contains the subgroup $R$ of rotations, so the
lift of a positive-energy projective unitary representation of
$L\T$ to an honest unitary representation of an extension
$\widetilde{L\T}$ of $L\T$ by a circle group $C$ restricts to a
representation of a semidirect product $R:E$, where $E$ is an
extension of $\Z \times \T$ by $C$ which splits over the identity
component. The irreducible positive-energy projective
representations of $L\T$ are classified by their restriction to
representations of $R \times C \times \T$  [30 \S 9.3]; an
irreducible representation of $\T$ is classified by its weight,
and the corresponding integer defined by $C$ is the level.
However, the identity component of $\widetilde{L\T}$wist has a
trivial subgroup of constant loops: its representation theory is
effectively weightless.

There is, however, an interesting basic representation
of $\widetilde{L\T}$wist; one construction, modelled on [34 \S
2], begins with the skew bilinear form defined on $L\R^{\Z/2\Z}$
by $$B(f_{0},f_{1}) = \frac{2}{\pi} \int_{0}^{1} f_{0}(x+\half)
f_{1}(x) dx \;.$$ The group $\D(2)$ of smooth 
orientation-preserving maps $g$ of $\R$ to itself satisfying
$g(x+\half) = g(x)+\half$ acts on this symplectic space, by $$g,f
\mapsto g^{\prime^{\half}}f \circ g.$$ The complexified space of
antiperiodic functions admits the decomposition
$$L\R^{\Z/2\Z}\otimes \C = A_{+} \oplus A_{-},$$ $A_{+}$ being
the subspace of functions on the circle which extend inside the
unit disk. There is a standard [34 \S 9.5] unitary representation
of $\widetilde{L\T}$wist on the symmetric algebra $\s (A_{+})$
associated to this polarization; the basis \medskip $$a_{n} = -
\pi i^{-\half} (-\half)^{n+\half} \Gamma(n+\half)^{-1}
\e((n+\half)x)$$ for the complexification satisfies
$$B(a_{n},a_{m}) = (2m+1)\delta_{n+m+1,0} \;.$$ The polarization is
defined by a nonstandard complex struction in which conjugation
acts by $\bar a_{n} = -ia_{-n-1}$, making $iB(\bar a,a)$ a
positive-definite Hermitian form on $A_{+}$. [This complex
structure differs from the standard one by a transformation which
is diagonal in the basis $a_{n}$; this operator is real but unbounded.] The
action of $\D(2)$ on $L\R^{\Z/2\Z}$ makes it reasonable to
interpret antiperiodic functions as sections of a bundle of 
half-densities on the circle; the complexification of this bundle
admits the nonvanishing flat section
$$(2\pi)^{\half}\e(x+\eighth)(dx)^{\half}
= (dZ)^{\half} \;,$$ where $Z = \e(x)$. It follows from Euler's
duplication formula that $$a_{n}(dx)^{\half} =
(2n+1)!!Z^{-n-1}(dZ)^{\half}$$ when $n$ is nonnegative.

The Lie algebra of $\D(2)$ now acts on $\s (A_{+})$ with
generators (cf.\ [20 \S 1.2, 37 \S 2])
$$\begin{array}{lllll}
L_{k} &=& \frac{1}{4} \sum_{n \in \Z}
a_{k-n-1}a_{n}\;\;&\mbox{if}\;\;& k \neq 0 \;,\\[10pt]
&=&\half \sum_{n\geq 0} a_{-n-1}a_{n} +
\frac{1}{16}\;\;&\mbox{if}\;\;&k = 0 \;.\end{array}$$ Convenient
polynomial generators $t_{n}$ for the algebra $\s (A_{+})$,
regarded as a ring of holomorphic functions on $A_{-}$, can be
defined by expanding an element $f$ of $A_{-}$ as $$\sum_{n\geq
0} t_{n}(f) Z^{-n-1}(dZ)^{\half} \;;$$ similar generators $T_{n}$,
constructed by writing this element as $$\sum_{n\geq 0} T_{n}(f)
a_{-n-1} \;,$$ satisfy the equation $$t_{n} = (2n+1)!! T_{n} \;.$$
The Virasoro generators (which are not derivations) act on these
elements so that
$$\begin{array}{lll}L_{k}T_{n} 
&=& (n-k+\half) T_{n-k}\;\;\;\mbox{if}\;\;\; n \geq k\\[10pt]
&=& 0 \;\;\;\mbox{ otherwise  .}\end{array}$$
On the other hand the group {\bf G} of invertible formal power
series (under composition) in $Z^{-1}$ acts on $A_{-}$,
interpreted as a free module over the ring of power series in
$Z^{-1}$; the Lie algebra of this group is spanned by vector
fields $$v_{k} = Z^{-k+1}d/dZ \;\;\;\mbox{ with }\;\;\; k \geq 0
,$$ which act on $\s (A_{+})$ (as derivations) such that
$$\begin{array}{lll}v_{k}t_{n} &=& (n-k+1)t_{n-k}\;\;\;\mbox{
when }\;\;\; n \geq k \\[10pt] &=& 0\;\;\;\mbox{  otherwise}
.\end{array}$$

We can thus identify $MU_{\C}^{*}(\pt)$ with $\s (A_{+})$
as a comodule over the Landweber-Novikov algebra, in a way which
makes it a Virasoro representation as well. [The complex
coefficients are only a technical convenience.] The resulting
bimodule defines a kind of Morita equivalence of the category of
cobordism comodules to the category of representations of
$\widetilde{L\T}$wist, given the (weak) monoidal structure
defined by the fusion product [35 \S 7] of representations.

\section*{Appendix II : the large phase space of deformations}\bigskip

Some recent advances [9] in quantum cohomology
seem to fit very naturally in the homotopy-theoretic framework sketched
above. This appendix has been added in July 1998 to the body of the published 
paper; I have used the opportunity to add references to recent
work on the construction of Gromov-Witten invariants discussed in \S 2.1.

\subsection*{II.1} If $H$ is a commutative Hopf algebra over a ring
$k$ and $B$ is a finitely-generated $\Z$-module, then the functor
defined on the category of commutative $k$-algebras by $$A
\mapsto B \otimes \Hom_{k-{\rm alg}}(H,A)$$ is represented by a Hopf
$k$-algebra which might be denoted $\otimes^{B} H$ : if $B$ is
free of rank $b$ then a choice of basis defines an isomorphism of
$\otimes^{B} H$ with the $b$-fold tensor product of copies of
$H$. A variant of this construction occurs in the theory of
vertex operator algebras; in that context $H$ is the
algebra of symmetric functions and $B$ is a positive even
lattice. 

This appendix suggests a conjectural interpretation for a family of
deformations of the quantum cohomology of a smooth algebraic
variety $V$ in terms of a similar construction, in which the role
of $B$ is played by the cohomology of $V$; in the basic example,
$H$ is an algebra of Schur $Q$-functions [28]. The parameter space for
this family is the `large phase space' of Witten [37 \S 3c], defined 
when topological gravity is coupled to quantum
cohomology; the more usual space of deformations of quantum
cohomology proper is then called the `small phase space'. In
impressionistic terms this large phase space is essentially a
tubular neighborhood of the moduli space of holomorphic maps from
a Riemann surface to $V$, inside the space of all smooth maps. A
two-dimensional quantum field theory is a kind of measure on such
a space of smooth (or perhaps continuous) maps, but it seems to
be reasonable to think of these measures as supported near the
(finite-dimensional) subspace of holomorphic maps. The resulting
hybrid structure can thus be interpreted as a homotopy-theoretic
family of deformations of a reasonably familiar kind of 
algebro-geometric object. 

This is a summary of work in progress; one way to paraphrase the
basic idea is that (even though quantum cohomology is not in any
very natural sense a functor), we might interpret it as a
cohomology theory taking values in the abelian category of
bicommutative Hopf algebras over $\Q$. The interesting examples
have further structure, and the point of this note is that some,
at least, of these extra structures seem to have natural
interpretations in a Hopf-algebraic context. 

I am indebted to Andy Baker and to Ezra Getzler for very helpful
discussions of this and related material. 

\subsection*{II.2} I will simplify notation here, writing $\M^{n}(V)$ for
$\M^{n}_{g,\alpha}(V)$ when indexing components is unnecessary. 
Gromov-Witten invariants in the context of algebraic geometry have now 
been defined rigorously by Behrend and Fantechi [2,3]; related results 
in the symplectic context have been announced by Li and Tian [23, cf.\ also 32]. 
We will be especially interested in the forgetful evaluation maps 
$$\epsilon^{n}(k) : \M^{n+k}(V) \rightarrow \M^{n}(V) \times V^{k} \; ;$$ 
the spaces of nonzero tangent vectors at marked points are principal
$\C^{\times}$-orbifold bundles over $\M^{n}(pt)$ which can be pulled back
over the domain of $\epsilon$, and $$v^{-k(d-1)}
\epsilon^{n}_{\T}(k) \in [\M^{n}(V) \times V^{k}_{+}, B\T^{k}_{+}
\wedge \MU_{\Lambda}]$$ will denote its cobordism class enriched
by the memory of the tangent bundles at the forgotten points.

Now let $$z = \sum t_{k,i} c^{k} z_{i} \in
H^{*}(B\T,H_{*}(V,\Q)) \; ,$$ where $\{z_{i}\}$ is a basis for
$H_{*}(V,\Q)$, be a homogeneous class of degree zero: we thus
interpret the coefficient $t_{k,i}$ to be an indeterminate of
degree $|z_{i}| - 2k$. In the language of physics, the elements
$z_{i}$ are the `primary fields' of a topological field theory
defined by the quantum cohomology of $V$, while $c^{k}z_{i}$ is
the $k$th `topological descendant' of $z_{i}$. Using this
terminology we can generalize the constructions of \S 2 above,
replacing the class $\phi^{n}_{g}(k)$ defined there with
$$\phi^{n}_{g}(k;z) = \sum _{\alpha \in H_{2}(V,\Z)} \Phi^{n}_{g,
\alpha *}(\epsilon^{n}_{\T,g,\alpha}(k) \cap \otimes^{k} z)
\otimes \alpha v_{(k)} \in H^{*}(\M^{n}_{g} \times V^{n},
\Lambda) \; ,$$ the cap product being (the $\Q$-linear extension
of) the Kronecker pairing $$H^{*}(B\T,H_{*}(V,\Z)) \otimes
H^{*}(V,H_{*}(B\T,\Z)) \rightarrow \Z \; .$$ The product with 
$z$ leaves the coproduct formula 2.1 essentially unchanged, and
the sum $$\tau^{n}_{g}(z) = v^{-d(n+g-1)} \sum_{k \geq 0}
\phi^{n}_{g}(k;z)$$ still defines a generalized topological field
theory, and thus a family of multiplications, which specializes
when $z$ is the fundamental class of $V$ to our previous construction. 
If, for example, $V$ is a single point, then $\tau^{0}(z)$
is a formal function from $H^{*}(B\T,\Q)$ to $H^{*}(\M,\Q)$,
which can alternately be interpreted as an element of the tensor
product of the symmetric algebra on $H_{*}(B\T,\Q)$ and the
cohomology of the moduli space of curves. On the other hand, the
composition $$B\T \rightarrow BU \rightarrow \MU$$ of the Thom
map for cobordism with the map induced by the inclusion of the
circle in the unitary group defines a canonical isomorphism
$$\s (H_{*}(B\T,\Q)) \rightarrow H_{*}(\MU,\Q) \; ;$$ the element
$\tau^{0}(z)$ thus defines a homomorphism from the homology of
the moduli space to $H_{*}(\MU,\Q) = MU_{*} \otimes \Q$. In this
case there is a unique primary field $z_{0} = 1$, and the
$t_{k}$'s defined by its topological descendents are the standard
generators of the Landweber-Novikov algebra.

\subsection*{II.3} At the opposite extreme we might suppose that $V$
is nontrivial and set $t_{i,k} = 0$ if $i,k > 0$; in this model
there are then nontrivial primary fields, but only the field
$z_{0} = 1$ has topological descendants. If we identify the
symmetric algebra on $H_{*}(B\T,\Q)$ with $MU_{*} \otimes \Q$ as
above, then $\tau^{n}(z)$ becomes the class in
$MU^{*}_{\Lambda}(V^{n})$ defined by the space of stable maps
from a curve marked with $n$ smooth points, together with an
indeterminate number of further distinct but unordered smooth
points which map to the subvariety $z$ of $V$. This yields a
deformation of the quantum multiplication on
$MU^{*}_{\Lambda}(V)$; its parameter space is the classical
`small phase space' of deformations, enlarged slightly (since we
permit nontrivial descendents of $z_{0}$) to yield a theory
interpretable in terms of cobordism rather than ordinary
cohomology. When $V$ is a point, for example, the deformation $z
\in MU^{*} \otimes \Q[v,v^{-1}]$ replaces the coupling constant
$q=q(v)$ of \S 2.3 with $q(zv)$.

In the general case (with no hypothesis that the topological
descendants of any primary fields vanish) we can interpret the
descendants of $z_{0} = 1$ as lying in the cobordism ring and
rewrite $\tau^{n}$ as an element of $$\s (H^{*}(V,H_{*}(B\T,\Q)))
\otimes H^{*}(\M^{n} \times V^{n}, \Lambda) \; ,$$ which can be
expressed in terms of reduced cohomology as $$\s (H_{*}(B\T,\Q))
\otimes \s ({\tilde H}^{*}(V,H_{*}(B\T,\Q))) \otimes H^{*}(\M^{n}
\times V^{n}, \Lambda) \; .$$ This is in turn isomorphic to the
symmetric $MU^{*}$-algebra $$\s_{MU_{\Q}} (MU^{*}_{\Q}(V))
\otimes_{MU_{\Q}} MU^{*}_{\Lambda}(\M^{n} \times V^{n}) \; ,$$
and we can think of $\tau^{n}(z)$ as a formal family of
deformations, parametrized by (the space underlying)
$H^{*}(V,H_{*}(B\T,\Q))$, of a generalized topological field
theory $$\tau^{n}(z) : \M^{n} \rightarrow {\bf F} (V^{n}_{+},
\MU_{\Lambda}) \; .$$

\subsection*{II.4}  Recently Eguchi [9, cf.\ also 13] and his coworkers have
studied an action of the Virasoro algebra on a large phase space
model for the quantum cohomology of $\C P(n)$, which generalizes
the Virasoro action on what can be
interpreted as the (large) quantum cohomology of a point. There
is reason to believe [29] that the latter Virasoro action can be
understood most naturally not in terms of complex cobordism but
instead in terms of cohomology with coefficients in a VOA-like
structure defined by the ring $\Delta$ of Schur $Q$-functions;
the resulting theory has good integrality properties, and its
Virasoro structure is a consequence of Hopf-algebraic properties
of $\Delta$. There is a natural (Kontsevich-Witten) genus $${\rm
kw}: MU_{*} \rightarrow \Delta[q_{1}^{-1}]$$ relating these
constructions, which is essentially an isomorphism over the
rationals, and it is natural to ask if this interpretation of
topological gravity can be extended to ecompass the quantum
cohomology of algebraic varieties. 

This is a subject for further research, but it is at least
reasonable in terms of the Hopf algebra structures. From that
point of view the map defined at the end of \S 3 identifies the rational
homology of $B\T$ with the primitives $P_{*}$ of $\Delta_{\Q}$,
which provides an interpretation of the large phase space of
deformations for the quantum cohomology of $V$ as the spectrum of 
$\Delta_{\Q}^{\otimes H(V,\Z)}$. Connected, graded Hopf algebras
over the rationals are primitively generated, and there is an
internal tensor product which sends the Hopf algebra $H_{0}$
(resp.\ $H_{1}$ with primitives $P_{0}$ (resp.\ $P_{1}$) to the
Hopf algebra with primitives $P_{0} \otimes P_{1}$; the Hopf
algebra of functions on the large phase space is just this tensor
product construction, applied to $\Delta_{\Q}$ and $\s (H^{*}(V))$. 
An internal tensor product on the category of bicommutative Hopf
algebras has become important recently in other parts of algebraic
topology, and
an integral version of the construction sketched above would seem
to be within reach. Questions related to Cartier duality and
self-adjointness need to be explored, but it seems likely that
these ideas will lead to a VOA-like
structure on this large quantum cohomology, at least in relatively simple
cases like $\C P(n)$.

\bibliographystyle{amsplain}

\begin{thebibliography}{99}

\bibitem [1]{1} P. Baum, W. Fulton, R. Mac Pherson, Riemann-Roch
for singular varieties, Publ. Math. IHES 45 (1975) 101-145

\bibitem [2]{2} K. Behrend, Gromov-Witten invariants in algebraic
geometry, Inv. Math 127 (1997) 601-617

\bibitem [3]{3} -----, B. Fantechi, The intrinsic normal
cone, Inv. Math, 128 (1997) 45-88 

\bibitem [4]{4} -----, Yu. Manin, Stacks of stable maps and
Gromov-Witten invariants, Duke Math J. 85 (1996) 1-60
 
\bibitem [5]{5} P. Berthelot et al, {\bf Th\'eorie des intersections et
theoreme de Riemann-Roch} [SGA 6] Lecture Notes in Mathematics 225, 
Springer (1971)

\bibitem [6]{6} M. Betz, J. Rade, Products and relations in
symplectic Floer homology, Stanford preprint (1995) 

\bibitem [7]{7} R.L. Cohen, J.D.S. Jones, G.B. Segal, Floer's
infinite dimensional Morse theory and homotopy theory, {\bf Floer
Memorial Volume}, Birkh\"auser, Progress in Mathematics 133 (1995)
297-326

\bibitem [8]{8} R. Dijkgraaf, E. Verlinde, H. Verlinde, Loop
equations and Virasoro constraints in nonperturbative
two-dimensional topological gravity, Nucl. Phys. B348 (1991)
435-456

\bibitem [9]{9} Eguchi T., Hori K., Xiong C-S., Quantum
cohomology and Virasoro algebra, Phys. Lett. B 402 (1997) 71-80

\bibitem [10]{10} Fukaya K., Morse homotopy and its quantization, {\bf 
Georgia Conference on Topology}, AMS Studies in Adv. Math. 1997, 409-440

\bibitem [11]{11} W. Fulton, R. Mac Pherson, {\bf Categorical framework
for the study of singular spaces}, Mem. AMS no. 243 (1981)

\bibitem [12]{12} -----, K. Pandharipande, Notes on stable maps and 
quantum cohomology, {\tt alg-geom}/9608011, Proc. of Symposia in Pure
Math., to appear

\bibitem [13]{13} E. Getzler, Intersection theory on $\M_{1,1}$ and
elliptic Gromov-Witten invariants, J. of the AMS 4 (1997) 973-998

\bibitem [14]{14} -----, M. Kapranov, Modular operads, Compositio
Math. 110 (1998 65-126

\bibitem [15]{15} H. Gillet, Intersection theory on algebraic
stacks and $Q$-varieties, J. Pure and Applied Algebra 34 (1984) 193-240

\bibitem [16]{16} J. Greenlees, J.P. May, {\bf Generalized Tate
cohomology}, Mem. AMS no. 543 (1995)

\bibitem [17]{17} T. Kimura, J. Stasheff, A.A. Voronov, The operad
structures of moduli spaces and string theory, Comm. Math. Phys.
171 (1995) 1-25

\bibitem [18]{18} F.F. Knudsen, The projectivity of the moduli
space of stable curves II, Math. Scand. 52 (1983) 161-199

\bibitem [19]{19} M. Kontsevich, Enumeration of rational curves
via torus actions, in {\bf The Moduli Space of Curves}, Birkh\"auser,
Progress in Mathematics 129 (1995) 335-368

\bibitem [20]{20} -----, Intersection theory on the
moduli space of curves and the matrix Airy function, Comm. Math.
Phys. 147 (1992) 1-23

\bibitem [21]{21} -----, Yu. Manin, Gromov-Witten
classes, quantum cohomology, and enumerative geometry, Comm.
Math. Phys. 164 (1994) 525-562

\bibitem [22]{22} -----, Yu. Manin, Quantum cohomology of
a product, Inventiones Math. 124 (1996) 313-339

\bibitem [23]{23} J. Li, G. Tian, Virtual moduli cycles and GW
invariants of general symplectic manifolds, {\tt alg-geom 9608032} 

\bibitem [24]{24} E. Looijenga, Intersection theory on
Deligne-Mumford compactifications [after Witten and Kontsevich],
Sem. Bourbaki no. 768 (1993)

\bibitem [25]{25} D. McDuff, D. Salamon, {\bf J-holomorphic curves and
quantum cohomology}, AMS University Lecture Series no. 6 (1995)

\bibitem [26]{26} Yu. Manin, Generating functions in algebraic
geometry and sums over trees, in {\bf The Moduli Space of Curves},
Birkh\"auser, Progress in Mathematics 129 (1995) 401-418

\bibitem [27]{27} J. Morava, Topological gravity and algebraic
topology, {\bf Proceedings of the Adams Symposium}, vol. II, LMS
Lecture Notes 176 (1992)

\bibitem [28]{28} -----,  Primitive Mumford classes,
Contemporary Math. 150 (1993) 291-302
 
\bibitem [29]{29} -----, Schur $Q$-functions and a
Kontsevich-Witten genus, in {\bf Homotopy theory via
algebraic geometry and representation theory}, ed. M.
Mahowald, S. Priddy, Contemporary Math, to appear
 
\bibitem [30]{30} A. Pressley, G. Segal, {\bf Loop groups}, Oxford
University Press (1986)

\bibitem [31]{31} A. Robinson, Derived tensor products in stable
homotopy theory, Topology 22 (1983) 1-18

\bibitem [32]{32} Y. Ruan, Higher genus symplectic invariants and
sigma models coupled with gravity, {\tt alg-geo} 9601005
 
\bibitem [33]{33} -----, G. Tian, A mathematical theory of
quantum cohomology, J. Diff. Geo. 42 (1995) 259-367

\bibitem [34]{34} G. Segal, Unitary representations of some 
infinite-dimensional groups, Comm. Math. Phys. 80 (1981) 301-342

\bibitem [35]{35} -----, Two-dimensional conformal field
theories and modular functors, in {\bf IXth International Congress on
Mathematical Physics}, ed. B. Simon, A. Truman, I.M. Davies,
Adam-Hilger (1988)

\bibitem [36]{36} R.W. Thomason, T. Trobaugh, Higher algebraic
$K$-theory of schemes and of derived categories, in {\bf The Grothendieck
Festschrift}, vol III, ed. P. Cartier et al, Progress in Mathematics 88
(1990) 247-436, Birkh\"auser

\bibitem [37]{37} E. Witten, Two-dimensional gravity and
intersection theory on moduli space, Surveys in Differential
Geometry 1 (1991) 243

\end{thebibliography}

\end{document}